\documentclass{amsart}

\usepackage{amssymb}

\newtheorem{theorem}{Theorem}[section]
\newtheorem{lemma}[theorem]{Lemma}
\newtheorem{proposition}[theorem]{Proposition}
\newtheorem{corollary}[theorem]{Corollary}

\theoremstyle{definition}
\newtheorem{definition}[theorem]{Definition}

\theoremstyle{remark}

\newcommand{\cL}{{\mathcal L}}
\newcommand{\cE}{{\mathcal E}}
\newcommand{\cP}{{\mathcal P}}
\newcommand{\cC}{{\mathcal C}}
\newcommand{\cN}{{\mathcal N}}

\newcommand{\bC}{{\mathbb C}}
\newcommand{\bR}{{\mathbb R}}

\numberwithin{equation}{section}

\begin{document}

\title{Lipschitz Extension Constants Equal Projection Constants}

\author{Marc A. Rieffel}
\address{Department of Mathematics\\
University of California\\
Berkeley, CA\ \ 94720-3840}
\curraddr{}
\email{rieffel@math.berkeley.edu}
\thanks{The research reported here was supported in part by 
National Science Foundation Grant DMS-0200591.}


\dedicatory{Dedicated to the memory of Gert Kj\ae g\aa rd Pedersen}


\begin{abstract}
For a Banach space $V$ we define its Lipschitz 
extension constant, $\cL\cE(V)$, 
to be the infimum of the constants $c$ such that for every metric 
space $(Z,\rho)$, every $X \subset Z$, 
and every $f: X \to V$, there is an extension, $g$, of $f$ to $Z$ 
such that $L(g) \le cL(f)$, 
where $L$ denotes the Lipschitz constant.  The basic theorem is 
that when $V$ is finite-dimensional we have $\cL\cE(V) = \cP\cC(V)$ 
where $\cP\cC(V)$ is the well-known projection constant of $V$.  
We obtain some direct consequences of this theorem, 
especially when $V =  M_n(\bC)$. We then apply known techniques
for calculating projection constants,
involving averaging of projections, 
to calculate $\cL\cE((M_n(\bC))^{sa})$.  
We also discuss what happens 
if we also require 
that $\|g\|_{\infty} = \|f\|_{\infty}$.
\end{abstract}

\maketitle

In my exploration of the relationship between 
vector bundles and Gromov--Hausdorff distance \cite{Rf} 
I need to be able to extend matrix-valued functions 
from a closed subset of a compact metric space to 
the whole metric space, with as little increase of the 
Lipschitz constant as possible.  
There is a substantial literature concerned with extending 
Lipschitz functions, 
but I have had difficulty finding there the facts which I need.  
The purpose of this largely expository paper is to describe and employ
a very strong relationship between the Lipschitz extension problem and 
what is referred 
to as the ``projection constant'' for finite-dimensional Banach spaces.  
This permits us to bring 
to bear on the Lipschitz extension problem the quite substantial 
literature concerning projection constants.  
This then provides the facts that I need, 
as well as other interesting facts.

In Section~\ref{sec1} we introduce what we call the 
Lipschitz extension constant, 
$\cL\cE(V)$, of a Banach space $V$. 
I have not found exactly this definition 
in the literature, although there are definitions very close to it.  
We also recall the well-known definition 
of the projection constant, $\cP\cC(V)$, 
of a Banach space $V$.  
The basic theorem is that if $V$ is finite-dimensional, then
$\cL\cE(V) = \cP\cC(V)$. I have not found this 
theorem stated in the
literature, probably because $\cL\cE(V)$ is not defined in the literature,
but I am told that this theorem is well-known to specialists on the geometry
of Banach spaces.

In Section~\ref{sec2} we give the proof that $\cL\cE(V) \le \cP\cC(V)$, 
while in Section~\ref{sec3} we 
give the proof that $\cP\cC(V) \le \cL\cE(V)$, thus proving the basic theorem.
We also show that restricting attention to compact metric space, or
to finite metric spaces, does not change this relation with the
projection constant.

In Section~\ref{sec4} we give some consequences of the basic theorem that
come directly from using facts about projection constants that are available
in the literature.  One of these consequences is that 
$\cL\cE(\bC) = 4/\pi$.  
Another of these consequences is the formula for $\cL\cE(M_n(\bC))$. Neither
of these consequences seems to have been recorded in the literature before.

However, what I specifically need for my exploration of 
vector bundles and Gromov--Hausdorff 
distance is $\cL\cE((M_n(\bC))^{sa})$, 
where $M_n(\bC)^{sa}$ denotes the Banach 
space of self-adjoint matrices in $M_n(\bC)$ with the operator norm. 
I have not seen
how to obtain this directly from facts stated in the literature. 
But in Sections 5 and 6
we discuss known techniques for dealing with projection constants, 
involving subspaces $V$
of $C(M)\ $for $M$ compact, and averaging of 
projections when there is a sufficiently
large group of isometries present. 
In Section 7 we then use these techniques to 
show that for any $n \geq 1$ we have
\[
\cL\cE(M_n(\bC)^{sa}) = 2n\Big(\frac{n}{n+1}\Big)^{n-1} \ \ -1   .
\]

Finally, in Section 8 we use radial retractions to discuss extending 
Lipschitz functions 
without increasing their supremum norm.

I am deeply indebted to Assaf Naor for patiently 
answering my occasional emailed
questions about this topic over the course of a number of months. 
He gave me important suggestions and brought to my attention 
important facts in the literature. I am equally deeply
indebted to William B. Johnson for his comments 
on the first version of this
paper and for patiently answering my 
subsequent emailed questions. He gave me further
important suggestions and brought to my 
attention further important facts in the
literature. I give some specific acknowledgments of their help at various 
points later in the paper. Much of this 
paper consists of little more than putting
together the items in the literature that they pointed out to me, 
and so this paper must
be considered largely expository. 

For the convenience of the reader I have 
included a number of the arguments that 
appear in the literature, and I have tried to formulate them
in a relatively constructive form. The audience that I
have had in mind when writing this paper consists of 
topologists and geometers who may
read my paper \cite{Rf} and would like to gain an 
understanding of the facts about
extending Lipschitz functions that I use there. More generally, 
this paper can be considered to
be an advertisement, detectable by MathSciNet and 
Google searches, that mathematicians 
who discover that
they need to extend vector-valued 
Lipschitz functions can find a body of facts
in the literature on the geometry of Banach 
spaces which may be quite useful to them. 

\section{The definitions and the main theorem}
\label{sec1}

Throughout this paper $(Z,\rho)$ will denote 
a metric space, and $X$ will denote a 
closed subset of $Z$ with the metric from $\rho$.  
Throughout this section we will assume 
that our Banach spaces are over $\bR$ unless the contrary is indicated, 
since for the Lipschitz extension problem it is irrelevant 
whether they are over $\bR$ or $\bC$.  
We will let $V$ denote a Banach space, often finite-dimensional.  
Let $f$ be a function from $X$ to $V$.  
Its Lipschitz constant, $L(f)$, is defined by
\[
L(f) = \sup\{\|f(x) - f(y)\|/\rho(x,y): x,y \in X,\ x \ne y\}.
\]
It can easily happen that $L(f) = +\infty$.  We define the Lipschitz 
constant of a function 
from $Z$ to $V$ similarly. In much of the discussion below we 
do not actually need to assume
that $X$ is closed, since if $L(f) < \infty$ then $f$ extends 
to the closure of $X$ with no
increase in Lipschitz constant. 
 
In general it is not possible to extend to $Z$ a function 
from $X$ to $V$ without increasing 
its Lipschitz constant.  I have not found the following 
definition in the literature, although 
there are definitions very close to it, such as $e(Y,Z)$ 
defined in section~1 of \cite{LN}.

\begin{definition}
\label{def1.1}
For a Banach space $V$ we let $\cL\cE(V)$ denote the 
infimum of the constants $c$ such that for 
any metric space $(Z,\rho)$ and any $X \subseteq Z$, 
and any function $f: X \to V$, 
there is an extension, $g$, of $f$ to $Z$ such that $L(g) \le cL(f)$.  
If no such constant $c$ 
exists, then we set $\cL\cE(V) = +\infty$.  We call $\cL\cE(V)$ 
the {\em Lipschitz extension constant} of $V$. We define $\cL\cE_c(V)$ 
much as above but
using only metric spaces $Z$ that are compact, and we
define $\cL\cE_f(V)$ much as above but using only metric spaces 
that are finite sets.
\end{definition}

Of course $\cL\cE_f(V) \leq \cL\cE_c(V) \leq \cL\cE(V)$. 

\begin{proposition}
\label{prop1.2}
Let $V$ be a finite-dimensional Banach space.  Then $\cL\cE(V) < \infty$.  
Furthermore, 
if $Z$ is compact then the infimum in the above 
definition is actually achieved, 
that is, for all $f$ there is an 
extension $g$ such that $L(g) \le (\cL\cE(V))L(f)$.
\end{proposition}

\begin{proof}
We use McShane's theorem \cite{MS} which states that every $\bR$-valued 
function on $X$ can be extended to $Z$ 
with no increase in the Lipschitz constant.  
We review its proof in the next section.  
Let $\{b_j\}$ be a basis for $V$ with $\|b_j\| = 1$ for all $j$, 
and let $\{\varphi_j\}$ be 
the dual basis.  Suppose that we are given $f: X \to V$.  
Set $f_j = \varphi_j \circ f$ for each $j$.  
For each $j$ choose, by McShane's theorem, 
an extension $g_j$ of $f_j$ to $Z$ such that 
$L(g_j) = L(f_j)$.  Define $g: Z \to V$ by $g(z) = \Sigma g_j(z)b_j$.  
Then $g$ is an extension of $f$ to $Z$, and
\[
L(g) \le \Sigma L(g_j) \le (\Sigma \|\varphi_j\|)L(f).
\]
Note that $\Sigma\|\varphi_j\|$ is independent of $Z, X$ and $f$.  
Thus $\cL\cE(V) < \infty$. (In fact, by Auerbach's 
lemma \cite {Wj} $\{b_j\}$ can be chosen
such that $\|\varphi_j\| = 1$ for all $j$, 
so that $L(g) \leq n L(f)$ where $n = \dim(V)$, and so $\cL\cE(V) \leq n$.)

Suppose now that $Z$ is compact. 
Choose some $c_1$ with $\cL\cE(V) < c_1 < \infty$.  
By the Arzela--Ascoli theorem one sees easily that, 
for a given $f: X \to V$ with $L(f) < \infty$, the set of 
its extensions $g: Z \to V$ for 
which $L(g) \le c_1L(f)$ forms a sup-norm compact subset 
of the functions from $Z$ to $V$.  
Furthermore $L$, as a function on the set of Lipschitz functions, 
can easily be verified to be 
lower semi-continuous for the sup-norm.  Thus there will be at 
least one extension $g$ for 
which $L(g)$ is minimal.  
It is then easily seen that $L(g) \le (\cL\cE(V))L(f)$.
\end{proof}

We now turn to the topic of projection constants \cite{Gr,IS,KT,Tm,KK}.

\begin{definition}
\label{def1.3}
For a Banach space $V$ we let $\cP\cC(V)$ denote the 
infimum of the constants $c$ such that for 
any Banach space $W$ into which $V$ is isometrically 
embedded there is a (linear) projection $P$ 
from $W$ onto $V$ such that $\|P\| \le c$.  
If no such constant $c$ exists then we set $\cP\cC(V) = +\infty$.
\end{definition}

\begin{proposition}
\label{prop1.4}
Let $V$ be a finite dimensional Banach space.  Then $\cP\cC(V) < \infty$.  
Furthermore, the infimum in 
the above definition is actually achieved, that is, 
for every $W$ containing $V$ there is a projection $P$ such 
that $\|P\| \le \cP\cC(V)$.
\end{proposition}

\begin{proof}
Let $\{b_j\}$ and $\{\varphi_j\}$ be as 
in the proof of Proposition~\ref{prop1.2}.  
Use the Hahn--Banach theorem to extend 
each $\varphi_j$ to a linear functional, 
${\tilde \varphi}_j$, 
on $W$ with $\|{\tilde \varphi}_j\| = \|\varphi_j\|$.  
Define $P: W \to V$ by 
$P(w) = \Sigma {\tilde \varphi}_j(w)b_j$.  
Then $P$ is a projection of $W$ onto $V$, 
and $\|P\| \le \Sigma \|\varphi_j\|$.  
Note that this bound is independent of $W$.  
Thus $\cP\cC(V) < \infty$. ( From Auerbach's lemma, 
stated above, we actually obtain
$\cP\cC(V) \leq n$.)

For a given $W$ it follows from the definition of $\cP\cC(V)$ 
that there is a sequence, $\{P_n\}$, 
of projections from $W$ onto $V$ such that 
$\|P_n\| \le \cP\cC(V) + 1/n$ for each $n$.  
Because $V$ is finite dimensional, 
the collection of operators $T$ from $W$ to $V$ for 
which $\|T\| \le k$ for some fixed constant $k$ is 
compact for the topology of pointwise 
convergence (by essentially the same proof 
as that of Alaoglu's theorem \cite{Rd2}, 
or by applying Alaoglu's theorem).  A limit, $P$, 
of the sequence $\{P_n\}$ is easily seen to be a 
projection of $W$ onto $V$ such that $\|P\| \le \cP\cC(V)$.
\end{proof}

I thank Assaf Naor for encouraging me to expect
that projection constants are relevant to the Lipschitz extension problem.

The basic theorem used in this paper is:

\begin{theorem}
\label{th1.5}
For any finite-dimensional Banach space $V$ we have 
\[\cL\cE(V) = \cP\cC(V) = \cL\cE_c(V) = \cL\cE_f(V).
\]
\end{theorem}

Thus we see that one benefit of 
introducing $\cL\cE_c(V)$ and $\cL\cE_f(V)$ is to see that if
one is working in a setting where one 
is only dealing with compact, or finite, metric spaces,
there is nevertheless no reduction in the Lipschitz extension constant. 
(They are also useful for technical purposes. See Theorem \ref{th5.1}.)

The above theorem is false in general for infinite-dimensional
Banach spaces. In fact, for every infinite-dimensional
separable Banach space $V$ we have $\cP\cC(V) = \infty$, for
(we paraphrase some lines on page 32 of \cite{BL}) 
Grothendieck showed in \cite{Grt} that every operator $T$
from $\ell^\infty(\Gamma)$ to a separable Banach space
is weakly compact. (Here $\Gamma$ is any discrete set.)
If $T$ is actually a projection onto a separable subspace of
$\ell^\infty(\Gamma)$, then from the Dunford-Pettis theorem
it follows that $T$ is actually compact, and so has 
finite-dimensional range. But every Banach space can be
isometrically embedded into some $\ell^\infty(\Gamma)$.
On the other hand, if $M$ is a compact and metrizable space then $C(M)$
is separable, and we will see that $\cL\cE_c(C(M)) = 1 $ 
whenever $M$ is compact. See also Theorem \ref{th5.3}
below.

In Section~\ref{sec8} we will consider the variation 
on $\cL\cE(V)$ in which we require of the 
extension $g$ not only that $L(g) \le cL(f)$ but also 
that $\|g\|_{\infty} = \|f\|_{\infty}$, 
where $\|\cdot\|_{\infty}$ denotes the supremum norm using the norm of $V$.

Since $\cL\cE_c(V) \leq \cL\cE(V)$, to prove Theorem \ref{th1.5} 
it suffices just 
to prove that $\cL\cE(V) \leq \cP\cC(V)$
and $\cP\cC(V) \leq \cL\cE_c(V) = \cL\cE_f(V)$.

\section{The proof that $\cL\cE(V) \le \cP\cC(V)$}
\label{sec2}

The basic extension theorem for functions with values 
in $\bR$ goes back to McShane in 1934 \cite{MS}.  
We sketch its proof.  It will be convenient here
to denote $\max\{r,s\}$ for $r,s \in \bR$ 
by $r\vee s$, and to use $\bigvee$ for the 
supremum of a bounded subset of $\bR$.  
We will also use these symbols for the max and supremum of a 
collection of $\bR$-valued functions on a set.

\begin{theorem}[McShane]
\label{th2.1}
Let $(Z,\rho)$ be a metric space, and let $X$ be a subset of $Z$.  
Let $f$ be an $\bR$-valued function on $X$.  
Then there is a (non-unique) extension, $g$, of $f$ to $Z$ such that
L(g) = L(f). If $\|f\|_\infty < \infty$, then we can arrange that also
$\|g\|_{\infty} = \|f\|_{\infty}$.

\end{theorem}

\begin{proof}
(See also theorem 1.5.6 of \cite{We}.) 
We can assume that $L(f) < \infty$.  
For each $x \in X$ define $h_x$ on $Z$ by $h_x(z) = f(x) - \rho(z,x)L(f)$.  
Then $L(h_x) = L(f)$ and $h_x|_X \le f$, while $h_x(x) = f(x)$.  
Let $g = \bigvee\{h_x: x \in X\}$.  If we pick some ``base-point" $x_0$, then
for any $z \in Z$ we have
\[
f(x) - f(x_0) \leq L(f)\rho(x, x_0) \leq L(f)(\rho(z, x) + \rho(z, x_0)),
\]
so that
\[
h_x(z) = f(x) - L(f)\rho(z, x) \leq f(x_0) + L(f)\rho(z, x_0).
\]
Thus $g(z) < \infty$. Furthermore,
$L(g) \le L(f)$ (see, e.g., proposition $1.5.5$ of \cite{We}), 
while $g|_X = f$, so that, in fact, $L(g) = L(f)$.
If $\|f\|_\infty < \infty$, 
then $h_x \leq f(x) \leq \|f\|_\infty$ for all $x$,
so that $g \leq \|f\|_\infty$. We can then replace $g$ by
 $g \vee(-\|f\|_{\infty})$ to obtain the desired extension of $f$
such that $\|g\|_\infty = \|f\|_\infty$.
\end{proof}

The above theorem fails already for complex-valued functions.  
Actually, this has nothing to 
do with the product of complex numbers, but rather involves just 
the fact that, as a Banach 
space, $\bC$ is $\bR^2$ with the Euclidean metric.  
A standard example (e.g., example $1.5.7$ 
of \cite{We}) consists of a $4$-point space 
$Z = \{\alpha,\beta,\gamma,\mu\}$, with $\alpha$, 
$\beta$ and $\gamma$ having distance $2$ 
from each other and distance $1$ to $\mu$.  
Let $X = \{\alpha,\beta,\gamma\}$, and let $f: X \to \bC$ 
have range exactly the $3$ cube-roots 
of $1$ (or the vertices of any equilateral triangle in $\bR^2$).  
Then it is easily checked that 
there is no extension $g$ of $f$ to $Z$ such that $L(g) = L(f)$.  
Basically this is due to the 
difference in curvature between the metric space $Z$ and the 
Euclidean space $\bR^2$ \cite{LPS}.  
(See also theorem $1.3$ of \cite{NPSS}).)  
Euclidean $\bR^2$ is ``flat'', while $Z$ is hyperbolic-like.

However, there do exist other Banach spaces 
to which Theorem~\ref{th2.1} generalizes, 
and these will be useful for the proof.  
Let $\Gamma$ be a discrete set, possibly infinite, even uncountable.  
We let $\ell^{\infty}(\Gamma)$ denote the Banach space of 
bounded real-valued functions on $\Gamma$ with the 
supremum norm.  The following is well-known.  
(See, e.g., lemma $1.1$ of \cite{BL}.)  

\begin{proposition}
\label{prop2.2}
Let $Z$, $\rho$, $X$ and $\Gamma$ be as above. 
Any function, $f$, from $X$ to $\ell^{\infty}(\Gamma)$ 
has an extension, $g$, to $Z$ such that $L(g) = L(f)$. In particular,
$\cL\cE(\ell^\infty(\Gamma)) = 1$. If $\|f\|_\infty < \infty$, then we can
arrange that also $\|g\|_{\infty} = \|f\|_{\infty}$.  
\end{proposition}

\begin{proof}
For each $\gamma \in \Gamma$ let $f_{\gamma}$ denote the $\bR$-valued 
function whose value 
at $x \in X$ is $f(x)$ evaluated at $\gamma$.  
Then $L(f_{\gamma}) \le L(f)$.
Let $g_{\gamma}$ be an extension of $f_{\gamma}$ to $Z$ 
as per Theorem \ref{th2.1}, 
so that $L(g_{\gamma}) = L(f_{\gamma})$.  
Define $g: Z \to \ell^{\infty}(\Gamma)$ by $g(z)(\gamma) = g_{\gamma}(z)$.  
Then it is easily verified that $g$ is the desired extension. 
If $\|f\|_\infty < \infty$,
then $\|f_{\gamma}\|_{\infty} \le \|f\|_{\infty}$ 
for each $\gamma$. Thus we can
choose the above $g_\gamma$'s such that 
$\|g_\gamma\|_\infty = \|f_\gamma\|_\infty$ for each
$\gamma$. The resulting $g$ will then satisfy 
$\|g\|_\infty = \|f\|_\infty$.
\end{proof}

{\sc Proof that $\cL\cE(V) \le \cP\cC(V)$.}  
Let $V$ be a finite-dimensional Banach space.  
Let $\Gamma$ be a subset of the unit ball 
of the dual space $V'$ such that for every $v \in V$ we have 
$\|v\| = \sup\{|\langle v,\gamma\rangle|: \gamma \in \Gamma\}$ .  
For example, $\Gamma$ can be all of the unit sphere, 
or a dense subset of the unit ball, 
or the set of extreme points of the unit ball.  
Then each element of $V$ can be viewed as 
a function on $\Gamma$ in the evident way, 
and this provides an isometric embedding of $V$ 
into $\ell^{\infty}(\Gamma)$.  By the definition of $\cP\cC(V)$ 
and Proposition \ref{prop1.4} there is a projection $P$ 
from $\ell^{\infty}(\Gamma)$ onto $V$ such that $\|P\| \le \cP\cC(V)$.

Let $Z$, $\rho$, $X$ be as earlier, and let $f: X \to V$.  
Through the above embedding we can 
view $f$ as having its values in $\ell^{\infty}(\Gamma)$ and 
this does not change $L(f)$ or $\|f\|_{\infty}$.  
Then according to Proposition \ref{prop2.2} we can find a 
function $h: Z \to \ell^{\infty}(\Gamma)$ 
such that $h|_X = f$ while $L(h) = L(f)$ 
(and $\|h\|_{\infty} = \|f\|_{\infty}$
if $\|f\|_\infty < \infty)$.  Set $g = P \circ h$.  
Then $g: Z \to V$ and $g|_X = f$, 
while $L(g) \le \|P\|L(h) \le \cP\cC(V)L(f)$.\hfill$ \qed$

\medskip
In principle, the above proof gives a constructive method for
producing extensions $g$ for which  $L(g) \le \cL\cE(V)L(f)$,
for a given $V$. We need only make one choice of an isometric
embedding of $V$ into an $\ell^{\infty}(\Gamma)$, and then
find one projection, $P: \ell^{\infty}(\Gamma) \to V$ with
$\|P\| \le \cP\cC(V)$. We can then proceed as in the
second paragraph of the above proof. The basic theorem
then shows that this gives $g$ with  $L(g) \le \cL\cE(V)L(f)$.
In fact, the above proof shows that $\cL\cE(V) \leq \|P\|$, and
so the basic theorem will imply the well-known fact
that $\|P\| = \cP\cC(V)$.

Notice, however, that the above proof does not 
give $\|g\|_{\infty} = \|f\|_{\infty}$ when $\|f\|_\infty < \infty$.  
We can only conclude that $\|g\|_{\infty} \le \cP\cC(V)\|f\|_{\infty}$.  
But we will see in Section~\ref{sec5} that we can arrange 
that $\|g\|_{\infty} = \|f\|_{\infty}$ 
at the cost of only knowing that $L(g) \le 2\cP\cC(V)L(f)$.

\section{The  proof that $\cP\cC(V) \le \cL\cE_c(V) = \cL\cE_f(V)$}
\label{sec3}

The proof of the inequality is a very 
minor reworking of the second proof of theorem $7.2$ of \cite{BL}, 
which is descended from the proof of theorem~$2$ of \cite{Lin}.  
I am indebted to Assaf Naor for telling 
me that the proof of theorem $7.2$ of \cite{BL} 
was what I needed here. We give an outline of the
proof. But first we remark that the fact that the inequality
$\cP\cC(V) \le \cL\cE(V)$ holds follows swiftly from 
corollary 1 to theorem 3 of \cite{Lin}. If $V$ is embedded 
isometrically in a Banach space $W$, then the identity
map from $V$ to itself will, by definition, have an extension,
$g$, to all of $W$ such that $L(g) \leq \cL\cE(V)$. Then
corollary 1 to theorem 3 of \cite{Lin}
implies that there is a projection, $P$, from
$W$ onto $V$ such that $\|P\| \leq \cL\cE(V)$.

Here is the outline of the proof that $\cP\cC(V) \le \cL\cE_c(V)$.
We must show that whenever $V$ is isometrically 
embedded in some Banach space $W$, 
then there is a projection, $P$, from $W$ onto $V$ 
such that $\|P\| \le \cL\cE_c(V)$.

Suppose first that $W$ is finite-dimensional, 
and that $W$ contains $V$ isometrically.  
For any $r \in \bR^+$ let $B^V(r)$ denote the closed 
ball about $0$ of radius $r$ in $V$, 
and similarly for $B^W(r)$.  Then $B^W(3)$ is a compact 
metric space which contains $B^V(3)$ 
as a closed subset.  Let $f$ be the identity map 
from $B^V(3)$ into $V$.  By the definition 
of $\cL\cE_c(V)$ there is a function, $g$, from $B^W(3)$ 
into $V$ such that $g(v) = v$ for 
$v \in B^V(3)$, and $L(g) \le (\cL\cE_c(V))L(f) = \cL\cE_c(V)$.

Next we smooth $g$ in the direction of $V$ by convolving it 
with a non-negative symmetric 
$C^{\infty}$ function 
on $V$ supported in the interior of $B^V(1)$ and of integral 1 
(for some choice of translation-invariant measure on $V$).  
The resulting function, $h$, when viewed as defined
on $B^W(2)$, satisfies $L(h) \le L(g)$ and $h(v) = v$ for
$v \in B^V(2)$. Furthermore, the derivatives of $h$ in $V$-directions
exist in the interior of $B^W(2)$ and are continuous there.  
Now choose a subspace $U$ of $W$ which is complementary to $V$, and choose a
sequence $\{\psi_n\}$ of $C^\infty$ functions on $U$ supported in
the interior of $B^U(1)$,  which forms 
an approximate $\delta$-function at $0$ in $U$. 
Convolve $h$ by each $\psi_n$ to
obtain a sequence $\{j_n\}$ of smooth $V$-valued 
functions, viewed as defined on $B^W(1)$, such
that $L(j_n) \le L(h) \le \cL\cE_c(V)$ for each $n$.  
Let $D_0j_n$ denote the total derivative of $j_n$ at $0 \in W$, 
so that $D_0j_n$ is a linear operator 
from $W$ to $V$ such that $\|D_0j_n\| \le \cL\cE_c(V)$.  
Because the $\psi_n$'s form 
an approximate $\delta$-function, one
finds that $(D_0j_n)(v)$ converges to $v$ for each $v \in V$.

The sequence $D_0j_n$ is contained in the ball of 
linear operators from $W$ to $V$ of norm no greater 
than $\cL\cE_c(V)$, which is compact.  Thus there is a 
subsequence which converges to an operator $P$.  
It is clear from the remarks just above that $Pv = v$ 
for any $v \in V$, and that 
$\|P\| \le \cL\cE_c(V)$.  
Thus $P$ is our desired projection from $W$ onto $V$.

Suppose now that $W$ is infinite-dimensional. 
For each finite-dimensional subspace $U$ of
$W$ which contains $V$ one can choose as above a projection, 
$P_U$, from $U$ onto $V$
such that $\|P_U\| \leq \cL\cE_c(V)$. An argument 
similar to the proof of Alaoglu's theorem, 
or of Proposition \ref{prop1.4} above, then yields a 
projection, $P$,  from $W$ onto $V$
such that $\|P\| \leq \cL\cE_c(V)$. 

We now turn to proving that $\cL\cE_c(V)  \leq \cL\cE_f(V)$, so that
they are equal. Let $(Z, \rho)$ be a compact metric space, and let
$X \subseteq Z$. Consider first the case in which $X$ is a finite
subset. Let $f: X \to V$. Let $S$ be a countable dense subset
of $Z$ containing $X$, and let $S$ be enumerated in such a way that
$X = \{s_j: 1 \leq j \leq n\}$. For each positive integer $k$ let
$S_k = \{s_j: 1 \leq j \leq k\}$. Then for each $k > n$ we can find, by
the definition of $\cL\cE_f(V)$, a function $g_k$ from $S_k$ to $V$
that extends $f$ and is such that $L(g_k) \leq \cL\cE_f(V)L(f)$.
It is easily seen that for each $k$ we have $\|g_k\|_\infty \leq r$
where
\[
r = \|f\|_\infty + (\mathrm{diameter}(Z))\cL\cE_f(V)L(f) .
\]
For each $k > n$ define $\hat{g}_k$ to be the extension of $g_k$
to all of $S$ which has value $0_V$ for each $s_j$ with $j > k$.
Because $V$ is finite-dimensional, the set of all functions from
$S$ into the $r$-ball in $V$ about $0_V$ is compact for the topology
of pointwise convergence, by Tychonoff's theorem. Since $S$ is
countable this topology is metrizable. Thus there is a subsequence,
say $\{\hat{g}_{k_m}\}$, which converges pointwise to a function, $\tilde{g}$,
from $S$ to $V$. It is easily seen that
$L(\tilde{g}) \leq \cL\cE_f(V)L(f)$. Then $\tilde{g}$ extends to a function,
$g$, from the completion, $Z$, of $S$, still with $L(g) \leq \cL\cE_f(V)L(f)$,
and $g$ is an extension of $f$.

Suppose now that $X$ is a subset of $Z$ which is not finite. Let $T$ 
be a countable dense subset of $X$, and enumerate its elements as a 
sequence
$\{t_j\}$. For each $n$ let $T_n = \{t_j: 1 \leq j \leq n\}$, and
let $f_n$ be the restriction of $f$ to $T_n$. Note that $L(f_n) \leq L(f)$.
As seen in the paragraph above, for each $n$ we can find an extension, $g_n$,
of $f_n$ to $Z$ such that $L(g_n) \leq \cL\cE_f(V)L(f)$. Thus the sequence
$\{g_n\}$ is equicontinuous. Also, $\|g_n\|_\infty \leq r$
for each $n$, where $r$ is as defined in the previous paragraph. Thus 
by the Arzela-Ascoli theorem the
sequence is totally bounded, and so has a subsequence which converges
uniformly to a $V$-valued function, $g$, on $Z$. (Actually, pointwise
convergence would suffice.) It is easily seen that $g$ is an 
extension of $f$ such that $L(g) \leq \cL\cE_f(V)L(f)$ as desired.
$\qed$

\section{Some Consequences}
\label{sec4}

In his recent book \cite{We} Weaver remarks just after example $1.5.7$ 
that it seems not to be known 
what is the smallest constant $c$ such that any $\bC$-valued function $f$ 
on a subset of a 
metric space can be extended to a function $g$ on $Z$ such 
that $L(g) \le cL(f)$.  That is, 
what is $\cL\cE(\bC)$?  But B.~Grunbaum showed in \cite{Gr} 
that $\cP\cC(\bR^2) = 4/\pi$ when $\bR^2$ 
is equipped with the Euclidean norm.  From Theorem~\ref{th1.5} we 
then immediately obtain:

\begin{corollary}
\label{cor4.1}
$\cL\cE(\bC) = 4/\pi$.
\end{corollary}

The reader will find it an easy and entertaining exercise to apply the
techniques that we describe in Sections 5 and 6 to give a proof
of this fact.

A non-obvious result of Kadec and Snobar (\cite{KS}; and 
see theorem $9.12$ of \cite{Tm} or theorem III.B.10
of \cite{Wj}) 
states that for any Banach space $V$ of real dimension $n$ 
we have $\cP\cC(V) \le \sqrt{n}$.  Thus:

\begin{corollary}
\label{cor4.2}
For any Banach space of real dimension $n$ we have $\cL\cE(V) \le \sqrt{n}$.
\end{corollary}

But the Kadec--Snobar theorem has been improved 
in \cite{KT}, so that if, for example, 
$n = 4$ then $\cP\cC(V) \le (2+3\sqrt{6})/5 < 2$, so 
that the same holds for $\cL\cE(V)$.  
We refer the reader to that paper for upper bounds for
other values of $n$.

My interest in this whole topic originated in my need for information 
about $\cL\cE((M_n(\bC)^{sa})$. But one
can first ask for $\cL\cE(M_n(\bC))$. 
I am indebted to William B. Johnson for telling
me that from theorem 5.6b of \cite{GL} one can deduce 
that $\cP\cC(M_n(\bC)) = (\cP\cC(\bC^n))^2$. 
Now Rutovitz \cite{Rt} showed that an inequality
for $\cP\cC(\bC^n)$ obtained by Grunbaum \cite{Gr}
is actually an equality. A proof along the lines that 
we will use in the next sections is given in
corollary III.B.16 of \cite{Wj}. One obtains:  
\[
\cP\cC(\bC^n) = n\int_{S^n} |z_1| d\lambda(z) = 
\Gamma(3/2)\Gamma(n+1)/\Gamma(n + 1/2) \geq (1/2)\sqrt{n\pi}  ,
\]
where $S_n$ is the unit sphere in $\bC^n$ 
and $\lambda$ is the rotationally invariant
measure of mass 1 on $S^n$. We thus obtain:

\begin{corollary}
\label{cor4.3}
For each $n \geq 2$ we have 
\[
\cL\cE(M_n(\bC)) = 
\Big(\Gamma(3/2)\Gamma(n+1)/\Gamma(n + 1/2)\Big)^2 \geq (\pi/4)n.
\]
\end{corollary}

It is interesting to note, in contrast, that if $D_n$ 
denotes the $*$-subalgebra of diagonal matrices 
in $M_n(\bC)$, then $D_n$ is isometric 
to $\ell_{\bC}^{\infty}(\Gamma_n)$ where $\Gamma_n$ is a set 
with $n$ points, and so $\cL\cE(D_n) = 4/\pi$, 
as follows easily from Corollary~\ref{cor4.1} and the 
proof of Proposition~\ref{prop2.2}.  Thus $\cL\cE(D_n)$ is 
independent of $n$, in contrast to $\cL\cE(M_n(\bC))$.

Finally, we now give an elementary argument which gives 
the Kadec--Snobar upper bound for $\cL\cE(M_n(\bC))$.  Specifically, 
for $A = M_n(\bC)$ and for each $m$ with $0 \le m \le n-1$, 
let $A_m$ be the linear subspace of matrices 
$\{t_{ij}\}$ such that $t_{ij} = 0$ unless $i-j = m \mod n$.  
Thus $A_0$ is our earlier $D_n$, the algebra of 
diagonal matrices.  For each $m$, multiplication by an appropriate 
permutation matrix carries $A_m$ 
isometrically onto $A_0$.  Thus $\cL\cE(A_m) = \cL\cE(A_0) \le \sqrt{2}$ 
by considering real and 
imaginary parts (or $= 4/\pi$ by the non-elementary Corollary~\ref{cor4.1}).  
But $A = \bigoplus_{m=0}^{n-1} A_m$ (with the precise 
relation between the norms being obscure), 
and from this we can at least obtain quickly that 
$\cL\cE(M_n(\bC)) \le n\sqrt{2}$ (or $\le 4n/\pi$).

\section{The usefulness of $V \subset C(M)$}
\label{sec5}

We now begin our discussion of techniques that will permit us 
to compute $\cL\cE((M_n(\bC))^{sa})$.  A standard fact about 
projection constants is that if $M$ is a compact metric space 
and if $V$ is a finite-dimensional subspace of $C(M)$, 
then $\cP\cC(V)$ is equal to the norm of any projection 
from $C(M)$ onto $V$ of minimal norm.  
(See, e.g., theorem~III.B.5 of \cite{Wj}.)  We will see in the 
next sections how this can be used.  I am much indebted to 
William~B.~Johnson for pointing out this path to me.

Here we will give a proof of this standard fact by using
$\cL\cE_c$ and then Theorem \ref{th1.5}. This
gives a somewhat constructive way of producing 
Lipschitz extensions with minimal increase of the Lipschitz 
constant. I have not seen the following theorem stated
in the literature.

\begin{theorem}
\label{th5.1}
Let $M$ be a compact space. For any compact
metric space $(Z, \rho)$, any subset $X \subseteq Z$, and
any function $f: X \rightarrow C(M)$ there is an
extension, $g$, of $f$ from $Z$ to $C(M)$
with $L(g) = L(f)$ and $\|g\|_\infty = \|f\|_\infty$.
In particular, $\cL\cE_c(C(M)) = 1$. 
\end{theorem}

\begin{proof}
We can assume that $L(f) < \infty$ and that $X$
is closed. Let $M_{dis}$ denote $M$ with the discrete
topology. We view $f$ as a function from $X$ to
$\ell^\infty(M_{dis})$. Then we construct an extension,
$\tilde g$, of $f$ from $Z$ to $\ell^\infty(M_{dis})$
in almost the same way as was used in the proof of
Proposition \ref{prop2.2} using Theorem \ref{th2.1}.
For each $x \in X$ we set
\[
h_x(z, m) = f(x, m) - L(f)\rho(z, x)
\]
for $z \in Z$ and $m \in M$. Define $\tilde g$, with
values in $\ell^\infty(M_{dis})$, by
${\tilde g}(z, m) = \bigvee_x \ h_x(z, m)$. Then
${\tilde g}(z, m) \leq \|f\|_\infty$ for all $z$ and
$m$, and $\tilde g$ is an extension of $f$ with
$L(\tilde g) = L(f)$.

What we need to show is
that $\tilde g(z) \in C(M)$ for all $z \in Z$. 
Since $X$ is compact, $\{f(x): x \in X\}$ is a compact
subset of $C(M)$, and so is equicontinuous by
the Arzela-Ascoli theorem. 
Then $\{h_x\}_{x \in X}$ is easily seen
to be an equicontinuous family of functions on
$Z \times M$. It follows that $\tilde g$, as a
function on $Z \times M$, is continuous, and
thus uniformly continuous. Consequently $\tilde g$,
as a function on $Z$, has values in $C(M)$, as
needed. We can now define $g$ by
$g(z) = \tilde g(z) \vee (-\|f\|_\infty)$ to
obtain the desired extension of $f$. (A similar idea
to the above proof is indicated in remark 3.3
of \cite{LR}.) 
\end{proof}

When we apply Theorem \ref{th5.1} in a way 
very similar to that in the second paragraph
of the proof that $\cL\cE(V) \leq \cP\cC(V)$, we quickly obtain:

\begin{theorem}
\label{th5.2}
Let $M$ be a compact space, and let $V$
be a finite-dimensional subspace of $C(M)$. Let
$P$ be a projection from $C(M)$ onto $V$. Then
$\cL\cE_c(V) \leq \|P\|$. Thus if $P$ is a projection
of minimal norm, then 
\[
\cP\cC(V) = \cL\cE(V) = \cL\cE_c(V) = \|P||.
\]
\end{theorem}

We remark that another approach to proving Theorem
\ref{th5.2} is to use peaking partitions of the identity,
as discussed in lemma 2.1 of \cite{MP}, to produce
in $C(M)$ isometric copies of $\ell^\infty(\Gamma)$
containing subspaces which converge to $V$ for
Banach-Mazur distance, where the $\Gamma$'s are
finite sets of increasing size. 

We also remark that from theorem 6b of \cite{Lin} and
its proof we quickly
obtain in much the same way as for Theorem \ref{th5.2}:

\begin{theorem}
\label{th5.3}
Let $(M, d)$ be any metric space, and let $C_u(M)$
be the Banach space of all bounded uniformly
continuous real-valued functions on $M$. Then
$\cL\cE(C_u(M)) \leq 37$.
\end{theorem}

See also theorem 1.6 of \cite{BL} and its proof.
 
\section{Averaging of projections}
\label{sec6}

Given enough symmetry and favorable circumstances, one can construct 
projections of minimal norm.  The following proposition is due to 
Rudin \cite{Rd}, but has antecedents for the circle group.  
See also theorem~III.B.13 of \cite{Wj}.  The proof is straightforward.

\begin{proposition}
\label{prop6.1}
Let $G$ be a compact group, and let $\alpha$ be a strongly-continuous 
representation of $G$ by isometries on a Banach space $W$.  Let $V$ be a 
subspace of $W$ which is $\alpha$-invariant, and suppose that there 
is a projection $Q$ from $W$ onto $V$.  
Define an operator $P$ from $W$ to $V$ by 
\[
Pw = \int_G \alpha_g(Q(\alpha_g^{-1}w))dg,
\]
where the Haar measure on $G$ gives $G$ measure $1$.  Then $P$ is a 
projection from $W$ onto $V$, and $\|P\| \le \|Q\|$.  Furthermore, $P$ 
is $\alpha$-invariant, in the 
sense that $\alpha_g \circ P \circ \alpha_g^{-1} = P$ 
for all $g \in G$.
\end{proposition}

Thus if we are in a situation in which the $\alpha$-invariant projection 
from $W$ onto $V$ is unique, then we will know that this projection
is a projection of 
minimal norm.  But this situation is easily described.  Let ${\tilde G}$ 
denote the set of equivalence classes of real irreducible representations 
of $G$.  They are all finite-dimensional.  For any $\gamma \in {\tilde G}$ 
a suitable multiple of its character, which we denote by $p_{\gamma}$, will 
be an idempotent in $L^1(G)$ for convolution.  (See \cite{Sm}, especially 
the appendix to III.5.)  For any strongly continuous
representation $\alpha$ of $G$ 
on a Banach space $W$ the operator $\alpha_{p_{\gamma}}$ that is the 
integrated form of $p_{\gamma}$ will be a projection from $W$ onto a 
subspace, $W_{\gamma}$, of $W$.  This subspace is the $\gamma$-isotypic 
component of $W$, in the sense that any irreducible $\alpha$-invariant 
subspace of $W$ on which the
representation of $G$ gives a representation isomorphic 
to $\gamma$ will be contained in $W_{\gamma}$.  The
$W_{\gamma}$'s are disjoint, 
and their algebraic direct sum is dense in $W$.  The kernel of the 
projection $\alpha_{p_{\gamma}}$ is the closure of the direct 
sum of all the other isotypic components.

Let $V$ be an $\alpha$-invariant subspace of $W$.  Then it too has 
isotypic components, $V_{\gamma}$, and $V_{\gamma} \subseteq W_{\gamma}$ 
for each $\gamma \in {\tilde G}$.

\begin{definition}
\label{def6.2}
We will say that an $\alpha$-invariant subspace $V$ of $W$ is 
$\alpha$-{\em full} if $V_{\gamma} = W_{\gamma}$ for each $\gamma$ 
for which $V_{\gamma} \ne \{0\}$.
\end{definition}

If the subspace $V$ is $\alpha$-full and if $P$ is an $\alpha$-invariant 
projection onto $V$, then the kernel of $P$ 
must contain all of the $W_{\gamma}$'s 
for which $V_{\gamma} \neq W_\gamma$, 
and the direct sum of these $W_{\gamma}$'s 
will be dense in the kernel of $P$.  Thus $P$ will be the unique 
$\alpha$-invariant projection onto $V$.  A few moments of contemplating 
actions by the one-element group, and then the general situation, 
shows that if $V$ is not $\alpha$-full then an $\alpha$-invariant 
projection onto $V$, if it exists, can not be unique.  Putting all of 
this discussion together, we obtain:

\begin{theorem}
\label{th6.3}
Let $\alpha$ be a strongly-continuous
representation of the compact group $G$ by 
isometries on a Banach space $W$.  Let $V$ be an $\alpha$-invariant 
subspace of $W$.  If a projection from $W$ onto $V$ exists, 
then an $\alpha$-invariant projection exists.  This $\alpha$-invariant 
projection is unique exactly if $V$ is $\alpha$-full.  If $V$ is 
$\alpha$-full, then this $\alpha$-invariant projection is a projection 
of minimal norm from $W$ onto $V$.
\end{theorem}

When this theorem is combined with Theorem~\ref{th5.2}, we obtain:

\begin{corollary}
\label{cor6.4}
Let $\alpha$ be a continuous action of $G$ on a compact space $M$, 
and let $\alpha$ also denote the corresponding 
representation of $G$ on $C(M)$.  
Let $V$ be a finite-dimensional $\alpha$-invariant subspace of $C(M)$ 
which is $\alpha$-full, and let $P$ be the (unique) $\alpha$-invariant 
projection of $C(M)$ onto $V$.  Then $\cL\cE(V) = \|P\|$.
\end{corollary}

We will make good use of this corollary in the next section.
From it we can also obtain a swift proof of Corollary \ref{cor4.1} 
as follows. Let $G = T = M$ where $T$ is the circle group, acting
on itself by translation, and so acting on $C(T)$. Embed $\bR^2$
into $C(T)$ by sending $(a,b) \in \bR^2$ to $f_{a,b} \in C(T)$
defined by
\[
f_{a,b}(t) = a\cos t \ + \ b\sin t  .
\]

\section{The calculation of $\cL\cE((M_n(\bC))^{sa})$}
\label{sec7}

Let $E = (M_n(\bC))^{sa}$, with $n \ge 2$ fixed throughout this section.  
We equip $E$ with the operator norm.  We will use the technique suggested 
by Corollary \ref{cor6.4} to calculate $\cL\cE(E)$.  
We let $\mathrm{tr}$ denote the 
(un-normalized) trace, and define a real-valued inner product on $E$ by 
$\langle a,b\rangle_E = \mathrm{tr}(ab)$.  Thus every linear functional on 
$E$ can be represented by an element of $E$.  We let $M$ denote the set 
of rank-one projections in $E$.  It is a compact subset of $E$ 
(and it corresponds to the set of extreme points of the state space of $E$).  
We define a linear mapping, 
$\varphi$, from $E$ into $C(M)$ by 
$\varphi_a(p) = \mathrm{tr}(ap)$ for $a \in A$ and $p \in M$. 
It is easily seen that $\varphi$ is isometric.  
(Here $C(M)$ denotes real-valued functions.)  Thus by Theorem~\ref{th5.2} 
we know that $\cL\cE(E)$ is equal to the norm of a projection of minimal norm 
from $C(M)$ onto $\varphi(E)$.

Let $\alpha$ denote the action of $G = SU(n)$ on $E$ by conjugation, 
and let $\beta$ denote the representation of $G$ on $C(M)$ coming from the 
action of $G$ on $M$ by conjugating rank-one projections.  
Then for any $u \in G$, $a \in E$ 
and $p \in M$ we have
\[
(\beta_u\varphi_a)(p) = \varphi_a(u^{-1}pu) = \mathrm{tr}(au^{-1}pu) = 
\mathrm{tr}(uau^{-1}p) = \varphi_{\alpha_u(a)}(p),
\]
so that $\varphi$ is $G$-equivariant, and $\varphi(E)$ is a $G$-invariant 
subspace of $C(M)$.

Now as a $G$-space $E$ decomposes into the direct sum of the subspace $E_0$ 
of scalar multiples of the identity and the subspace $E_1$ of elements of 
trace $0$.  The latter is isomorphic as $G$-space, via multiplication 
by $i = \sqrt{-1}$, to $su(n)$, the Lie algebra of $G$, with the 
representation of $G$ on $su(n)$ being the adjoint 
representation (which is irreducible since $su(n)$ is a simple Lie algebra).  
Thus the representation 
of $G$ on $\varphi(E)$ is isomorphic to the direct sum of 
the trivial representation and the adjoint representation.  
Let $p_0 = \begin{pmatrix} 1 & 0 \\ 0 & 0_{n-1} \end{pmatrix}$.  
The stability subgroup, $H$, of $p_0$ consists of the $u$'s in $G$ 
of the form $u = \begin{pmatrix} \lambda & 0 \\ 0 & u' \end{pmatrix}$ 
with $u' \in U(n-1)$ and ${\bar \lambda} = \mathrm{det}(u')$.  
Thus we can identify $H$ with $U(n-1)$.  The action of $G$ on $M$ is 
transitive, and so we can identify $M$ with $G/H$.

The representation of $G$ on $C(M)$ 
can be viewed as the representation of $G$ 
obtained by inducing to $G$ the trivial representation of $H$ in the way 
described in section~III.6 of \cite{BD}.  The Frobenius reciprocity theorem, 
in the form given in proposition III.6.2 of \cite{BD}, 
shows that the multiplicity in $C(M)$ of an irreducible representation 
of $G$ will equal the multiplicity of the trivial representation of $H$ 
in the restriction to $H$ of that irreducible representation of $G$.  
When the representation of $G$ on $E_1$ is restricted to $H$, 
its subspace of invariant 
vectors consists of the scalar multiples of 
$\begin{pmatrix} n-1 & 0 \\ 0 & -I_{n-1} \end{pmatrix}$, 
and so is $1$-dimensional.  We then see that the adjoint representation 
and the trivial representation of $G$ in $C(M)$ 
each occur with multiplicity $1$.  
Thus $\varphi(E)$ is the direct sum of two isotypic components of the 
action of $G$ on $C(M)$.  Consequently $\varphi(E)$ is a full $G$-invariant 
subspace of $C(M)$.  
From Corollary~\ref{cor6.4} we see that the $G$-invariant 
projection from $C(M)$ onto $\varphi(E)$ is unique, and that its norm 
equals $\cL\cE(E)$.  We denote this projection by $P$.  We now need to 
find a convenient expression for $P$ that will permit us to calculate $\|P\|$.

From the discussion of multiplicity given above it is clear that the vector 
space of $G$-invariant operators from $C(M)$ onto $\varphi(E)$ 
is $2$-dimensional.  The projection, $P_0$, of $C(M)$ onto the 
constant functions is given by
\[
(P_0f)(p) = \int_M f(q)dq,
\]
where we use the $G$-invariant measure on $M$ which gives $M$ measure $1$.  
It is easily seen that the operator $T$ defined by
\[
(Tf)(p) = \int_M \langle p,q\rangle_Ef(q)dq = 
\langle p,\int_M qf(q)dq\rangle_E
\]
has values in $\varphi(E)$ and is $G$-invariant.  
Thus $P$ must be a linear combination of $P_0$ and $T$.  
So we will seek $P$ in the form $(Pf)(p) = \int_M K(p,q)f(q)dq$ 
where $K$ has the form
\[
K(p,q) = \mu + \nu\langle p,q\rangle_E,
\]
where $\mu$ and $\nu$ are constants to be determined.  
We must of course have $P1 = 1$.  By considering $\int_M \varphi_a(p)dp$ and 
its $G$-invariance, it is easy to see 
that $\int_M qdq$ is $\alpha$-invariant.  
Thus it must be a multiple of $I_n$.  On taking the trace we find that 
$\int qdq = n^{-1}I_n$.  From $P1=1$ it then follows that $\mu = 1-n^{-1}\nu$.

To obtain a second equation for $\mu$ and $\nu$ we use the reproducing-kernel 
property of $K$.  For $f,g \in C(M)$ we let $\langle f,g\rangle_M$ denote 
their usual inner product for $L^2(M)$.  For each $p \in M$ define 
$K_p$ by $K_p(q) = K(p,q)$.  Then the definition of $P$ can be 
rewritten as $(Pf)(p) = \langle K_p,f\rangle_M$.  Let $\{e_j\}_{j=1}^{n^2}$ 
be an orthonormal basis for $\varphi(E) \subseteq L^2(M)$.  
Then $K_p = \sum \langle K_p,e_j\rangle_Me_j$ for each $p \in M$.  But
\[
\langle K_p,e_j\rangle_M = (Pe_j)(p) = e_j(p).
\]
Thus
\[
K(p,q) = \sum e_j(p)e_j(q).
\]
Consequently
\[
\int_M K(q,q)dq = \sum \langle e_j,e_j\rangle_M = n^2.
\]
But from the definition of $K$ it is clear that $K(q,q)$ has the constant 
value $\mu+\nu$, so that $\mu+\nu = n^2$.  From this and the equation 
$\mu = 1 - n^{-1}\nu$ obtained earlier 
we find that $\mu = -n$ and $\nu = n(n+1)$.  

We must now determine $\|P\|$, 
where we now revert to the sup-norm on $C(M)$.  
Now $P$ is $G$-invariant and $G$ acts 
transitively on $M$, so it suffices to determine the norm of the linear 
functional $f \mapsto (Pf)(p_0)$ on $C(M)$.  In view of the form of $P$ 
this will be given by
\[
\int_M |-n + n(n+1)\langle p_0,q\rangle_E| dq.
\]
To evaluate this integral we use the following specialization to $\bC$ of 
lemma~$3'$ii of \cite{Rg2} (or lemma~$4.4$ of \cite{Rg1}):

\begin{lemma}
\label{lem7.1}
For any continuous function $h$ on the interval $[-1,1]$ we have
\[
\int_M h(2\langle p,q\rangle_E - 1)dq = 
(n-1)2^{1-n} \int_{-1}^1 h(t)(1-t)^{n-2}dt.
\]
\end{lemma}

But
\[
-n + n(n+1)\langle p,q\rangle_E = 
2^{-1}n((n-1) + (n+1)(2\langle p,q\rangle_E -1)),
\]
and so from Lemma~\ref{lem7.1} we see that
\[
\|P\| = 2^{-n}n(n-1) \int_{-1}^1 |(n-1) + (n+1)t|(1-t)^{n-2}dt.
\]
Straightforward calculation of this integral yields $\|P\| = 
2n\left( \frac {n}{n+1}\right)^{n-1} - 1$.  We thus obtain:

\begin{theorem}
\label{th7.2}
$\cL\cE((M_n(\bC))^{sa}) = 2n\left( \frac {n}{n+1}\right)^{n-1} -1$.
\end{theorem}

We can rewrite this formula as
\[
2(n+2+n^{-1})\left( 1 - \frac {1}{n+1}\right)^{n+1} - 1
\]
and notice that $\left( 1 - \frac {1}{n+1} \right)^{n+1}$ converges 
to $e^{-1}$ as $n \to \infty$.  If we rewrite our formula in the form
\[
\cL\cE((M_n(\bC))^{sa}) = n\omega(n),
\]
we then find that $\omega(n)$ converges to $2e^{-1}$ as $n \to \infty$.  
We also note that our formula gives the correct answer for $n=1$.

\section{Preservation of the supremum norm}
\label{sec8} 

We mentioned earlier that the extensions of Lipschitz 
functions that we have discussed so far do not 
always preserve the supremum norms of the functions. 
We will now show that 
the norm can be preserved at the cost of no 
more than doubling the Lipschitz constant. In
Theorems \ref{th2.1} and \ref{th5.1} we were able to use
the lattice structure of $\ell^\infty(\Gamma)$ and $C(M)$
and its relation to the norm in order to arrange 
preservation of the norm, but this technique is not 
generally available.

The tool which we use is that of radial 
retractions (which are also called ``radial projections'').  
Let $V$ be any normed vector space over $\bR$ or $\bC$.  
For any $r \in \bR$ with $r > 0$ define the 
radial retraction $\Pi_r$, a non-linear map from $V$ into itself, by
\[
\Pi_r(v) = \begin{cases}
v &\mbox{if $\|v\| \le r$} \\
rv/\|v\| &\mbox{if $\|v\| \ge r$.}
\end{cases}
\]
Let $L(\Pi_r)$ denote the Lipschitz constant of $\Pi_r$.  
The first assertion of the following proposition 
is basically known.  See \cite{DW} and \cite{Th}.  
I have not seen the next two assertions stated in 
the literature, though they can be obtained via theorem~$2$ of \cite{Th}. 

\begin{proposition}
\label{prop8.1}
For any normed vector space $V$ we have $L(\Pi_r) \le 2$.  
If $V = C(X)$ for any compact space $X$ 
containing at least two points, then $L(\Pi_r) = 2$.  
If $V$ is a $C^*$-algebra of dimension at 
least $2$, then $L(\Pi_r) = 2$.
\end{proposition}

\begin{proof}
We have $\Pi_r(v) = r\Pi_1(r^{-1}v)$ for $v \in V$, 
and so $L(\Pi_r) = L(\Pi_1)$.  Thus it suffices 
to prove the first assertion for $\Pi = \Pi_1$, 
which we now do.  If $\|v\| \le 1$ and $\|w\| \le 1$ 
then clearly $\|\Pi(v) - \Pi(w)\| = \|v-w\|$.  
Suppose that $\|v\| \ge 1$ while $\|w\| \leq 1$.  Then
\[
\|v/\|v\| - w\| \le \|v/\|v\| - v\| + \|v-w\| = \|v\| - 1 + \|v-w\|.
\]
But $\|v\|-1 \le \|v-w\| + \|w\| - 1 \le \|v-w\|$, 
so that $\|\Pi(v)-\Pi(w)\| \le 2\|v-w\|$.  Finally, 
suppose that $\|v\| \ge 1$ and $\|w\| \ge 1$.  
By symmetry we can assume that $\|w\| \le \|v\|$.  
Then $\|v/\|w\|\| \ge 1$ and $\|w/\|w\|\| \le 1$, 
so that we can apply the previous case to obtain
\[
\|(v/\|w\|)/\|v/\|w\|\| - w/\|w\|\| \le 2\|v/\|w\| - w/\|w\|\|.
\]
Upon simplifying, we obtain the desired inequality.

Suppose now that $V = C_{\bR}(X)$ for $X$ compact, 
and that $x$ and $y$ are distinct points of $X$.  
Choose $g \in V$ such that $-1 \le g \le 1$, $g(x) = 1$ 
and $g(y) = -1$.  For any $\epsilon > 0$ 
set $f = g + \epsilon$.  Thus $\|g\|_{\infty} = 1$, 
while $\|f\|_{\infty} = 1 + \epsilon$ and 
$\|f-g\|_{\infty} = \epsilon$.  But
\[
f/\|f\| - g = (g+\epsilon)/(1+\epsilon) - g = \epsilon(1-g)/(1+\epsilon).
\]
Now $\|1-g\|_{\infty} = 2$, and so 
$\|f/\|f\| - g\|_{\infty} = 2\epsilon/(1+\epsilon)$.  If $k$ is a 
constant such that
\[
\|f/\|f\| - g\|_{\infty} = \|\Pi(f) - \Pi(g)\|_{\infty} \le k\|f-g\|_{\infty},
\]
then the calculation above shows 
that $2\epsilon/(1+\epsilon) \le k\epsilon$.  
On letting $\epsilon$ go 
toward $0$, we see that $2 \le k$.  (We remark that the 
above argument works for any linear subspace 
$V$ of $C(X)$, or of $\ell^{\infty}(\Gamma)$, which 
contains the constant functions and at least one 
non-constant function, i.e. for order-unit spaces.)

Now any unital $C^*$-algebra, commutative or not, will, 
if it has dimension at least $2$, contain a 
commutative $C^*$-algebra of dimension at least $2$ and 
thus also a $C_{\bR}(X)$ for an $X$ with at 
least $2$ points.  Since $\Pi$ will carry linear subspaces 
into themselves, we can apply what we have 
found for $C_{\bR}(X)$ to conclude that $L(\Pi) = 2$ for 
unital $C^*$-algebras.  A bit more arguing 
deals with non-unital $C^*$-algebras.
\end{proof}

We remark that for a Hilbert space we have $L(\Pi) = 1$.

Suppose now that we have a metric space $(Z,\rho)$, 
a subset $X$, and a bounded function $f$ from $X$ 
into a Banach space $V$.  If $g$ is an extension of $f$ to $Z$, 
and if $r = \|f\|_{\infty}$, 
then $h = \Pi_r \circ g$ will be an extension of $f$ to $Z$ 
such that $\|h\|_{\infty} = \|f\|_{\infty}$ 
and $L(h) \le L(\Pi_r)L(g) \le 2L(g)$.

We can formalize the situation with the following definition.

\begin{definition}
\label{def5.2}
For a Banach space $V$ we let $\cL\cN\cE(V)$ 
denote the infimum of the constants 
$c$ such that for any  metric space $(Z,\rho)$ and 
any $X \subseteq Z$, and any bounded function 
$f: X \to V$, there is an extension, $g$, of $f$ to $Z$ such that 
$L(g) \le cL(f)$ and $\|g\|_{\infty} = \|f\|_{\infty}$. We
define $\cL\cN\cE_c(V)$ and $\cL\cN\cE_f(V)$ similarly.
\end{definition}

Then our discussion above gives:

\begin{proposition}
\label{prop5.3}
For any Banach space $V$ we have 
$\cL\cN\cE(V) \le 2\cL\cE(V)$.
\end{proposition}

Note that Proposition~\ref{prop2.2} says 
that $\cL\cN\cE(\ell_{\bR}^{\infty}(\Gamma)) = 1$ and Theorem \ref{th5.1} 
says that $\cL\cN\cE_c(C(M)) = 1$.  
It would be interesting to know the exact value
for $\cL\cN\cE(V)$ for various choices of $V$, 
expecially for $M_n(\bC)$ and $(M_n(\bC))^{sa}$. I have
not found a $V$ for which I could prove that $\cL\cN\cE(V) \neq \cL\cE(V)$.
We could look for other retractions of $V$ onto its unit ball. Thus
we seek information about $\cL{\mathcal R}(V)$, where $\cL{\mathcal R}(V)$
is the infimum of the constants $c$ such that there is a retraction, $R$,
from $V$ onto its unit ball with $L(R) \leq c$. We are thus
looking for Lipschitz extensions to the metric space $V$ of the function
consisting of the inclusion of the unit ball of $V$ into $V$. 
Proposition~\ref{prop8.1} shows that we always have 
$\cL{\mathcal R}(V) \leq 2$.

Assume that $V$ is finite-dimensional. Then arguments very similar 
to those towards the end of Section 3 show 
that $\cL\cN\cE_c(V) = \cL\cN\cE_e(V)$;
and somewhat similar argunents, involving the directed sets of
all compact (or finite) subsets of $Z$ and of $X$ and Tychonoff's
theorem, show that $\cL\cN\cE(V) = \cL\cN\cE_c(V)$.

\def\cprime{$'$}
\providecommand{\bysame}{\leavevmode\hbox to3em{\hrulefill}\thinspace}
\providecommand{\MR}{\relax\ifhmode\unskip\space\fi MR }
\providecommand{\MRhref}[2]{%
  \href{http://www.ams.org/mathscinet-getitem?mr=#1}{#2}
}

\def\cprime{$'$}
\providecommand{\bysame}{\leavevmode\hbox to3em{\hrulefill}\thinspace}
\providecommand{\MR}{\relax\ifhmode\unskip\space\fi MR }
\providecommand{\MRhref}[2]{%
  \href{http://www.ams.org/mathscinet-getitem?mr=#1}{#2}
}
\providecommand{\href}[2]{#2}

\end{document}